\documentclass[12pt]{article}%16.08.2015
\usepackage{amsmath}
\usepackage{amssymb}
\usepackage{dsfont}
\usepackage{cite}
\newsymbol \blackbox 1004
\newcommand{\eh}{\hfill}\newlength{\sperr}

\newenvironment{proof}{{\settowidth{\sperr}{\bf\rm
Proof}%
\par\addvspace{0.3cm}\noindent\parbox[t]{1.3\sperr}
{\bf\rm P\eh r\eh o\eh o\eh f\eh }%
}}{\nopagebreak\mbox{}
$\blackbox$\par\addvspace{0.3cm}}

\def\nn{\nonumber}
\def\a{\alpha}
\def\b{\beta}

\def\ka{\kappa}

\def\vr{\varrho}
\def\Lam{\Lambda}
\def\s{\sigma}
\def\la{\lambda}
\def\om{\omega}

\def\t{\theta}

\def\vp{\varphi}
\def\vt{\vartheta}
\def\ve{\varepsilon}
\def\wh{\widehat}
\def\wt{\widetilde}
\def\ov{\overline}

\def\BC{{\mathbb C}}
\def\BR{{\mathbb R}}
\def\BN{{\mathbb N}}

\def\cla{{\mathcal A}}
\def\clb{{\mathcal B}}
\def\clc{{\mathcal C}}

\def\clg{{\mathcal G}}
\def\clh{{\mathcal H}}

\def\clt{\mathcal{T}}

\def\im{{\rm Im\ }}

\def\diag{\mathrm{diag}}
\newcommand{\E}{\mathrm{e}}
\newcommand{\I}{\mathrm{i}}

\newtheorem{Pa}{Paper}[section]
\newtheorem{Tm}[Pa]{{\bf Theorem}}
\newtheorem{La}[Pa]{{\bf Lemma}}

\newtheorem{Rk}[Pa]{{\bf Remark}}

\newtheorem{Dn}[Pa]{{\bf Definition}}

\newtheorem{Pn}[Pa]{{\bf Proposition}}

\newenvironment{dedication}
        {\vspace{1ex}\begin{quotation}\begin{center}\begin{em}}
        {\par\end{em}\end{center}\end{quotation}}

\title{Inverse problems for self-adjoint Dirac systems: explicit solutions and stability of the procedure}

\author{A.L. Sakhnovich}

\date{}
\begin{document}
\maketitle

\begin{dedication}
To the memory of Leiba Rodman,
a wonderful mathematician and an admirable person.
\end{dedication}

\begin{abstract} A procedure to recover explicitly self-adjoint matrix Dirac systems on semi-axis
(with both discrete and continuous components of spectrum) from rational Weyl functions
is considered. Its stability is proved. GBDT version of B\"acklund-Darboux transformation
and various important results on Riccati equations are used for this purpose.
\end{abstract}



\section{Introduction}
\setcounter{equation}{0}
Self-adjoint Dirac system has the form
\begin{align} &       \label{1.1}
\frac{d}{dx}y(x, z )=\I (z j+jV(x))y(x,
z ), \qquad
x \geq 0,
\end{align} 
where
\begin{align} &   \label{1.2}
j = \left[
\begin{array}{cc}
I_{m_1} & 0 \\ 0 & -I_{m_2}
\end{array}
\right], \hspace{1em} V= \left[\begin{array}{cc}
0&v\\v^{*}&0\end{array}\right],  \quad m_1+m_2=:m,
 \end{align} 
$I_{m_k}$ is the $m_k \times m_k$ identity
matrix and $v(x)$ is an $m_1 \times m_2$ matrix function.
We assume that the {\it potential} $v$ is locally summable (i.e., summable on all the finite intervals $[0,l]$).

Inverse spectral problem for general-type self-adjoint Dirac system, and closely related problem to recover Dirac system
from its Weyl-Titchmarsh (Weyl) function, had been actively studied since 1950s \cite{Kre1, LS} and many interesting results were 
published last years (see, e.g.,  \cite{AGKLS2, BrEKT, ClGe, CG2,  EGNST, Hryn, MP, SaSaR}
and various references therein). Further, speaking about inverse 
spectral problems, we mean (in particular) the recovery of systems from their Weyl functions.  Inverse spectral problems (in the mentioned above sence) are also solved  \cite{FKRS2012, FKRS2014, SaA021, SaA14, SaSaR} for 
general-type skew-self-adjoint and discrete Dirac systems.

Procedures to solve these inverse problems are nonlinear and usually unstable. However, stability plays an essential role in theory and applications,
and several special cases where such stability could be proved are important. In particular, we could mention the paper \cite{MORos} on the evolution
Schr\"odinger equation and the paper \cite{Hryn}, where stability was proved for  a class of scalar ($m_1=m_2=1$) Dirac systems on interval  with discrete and
$d$-separated spectral data. Here, we consider the case of explicit solutions of inverse problems (i.e., the case of rational Weyl functions), and one can apply  procedures, which are
different from the procedures for the general-type case. 

Thus, we prove stability of solving inverse problems for matrix Dirac systems on semi-axis
with both discrete and continuous components of spectrum.
Riccati equations play in essential role in the explicit solving of inverse problems,
and so the classical works on Riccati equations by Leiba Rodman and coauthors are actively used in this paper.

Explicit solutions of inverse spectral problems can be obtained either by applying procedures for general-type systems to some
specific spectral data (e.g., to rational Weyl or scattering functions) or by using several specific for explicit solutions procedures.
The first (general-type) approach was used, for instance, in \cite{AlGo, FKS1, ASAK} and \cite[Sect. 6]{AGKLS2}. The second
approach includes Crum-Krein method \cite{Cr, KrChr}, commutation methods \cite{D, Ge, GeT, KoSaTe} and  some versions of B\"acklund-Darboux transformation.
Here we use our GBDT version of the B\"acklund-Darboux transformation (see \cite{SaA1, SaAmmnp10, SaSaR} and references therein), see also
\cite{C1, Ci, MS, ZM} and references therein on various versions of  B\"acklund-Darboux transformations.

In the next section, Preliminaries, we present some basic notions from system theory and formulate several results on Weyl functions.
We also present GBDT  procedure to solve inverse problem for systems \eqref{1.1} (more precisely, to recover self-adjoint Dirac systems from Weyl functions).
Section~\ref{Sta} is dedicated to the proof of stability of this procedure.

As usual, $\BR$ stands for the real axis, $\BC$ stands for the complex plane, $\BC_+$ is the open upper half-plane
$\{z: \,\Im(z)>0\}$, $\BC_-$ is the open lower half-plane $\{ z: \, \Im (z)<0\}$, and the notation diag$\{d_1, ...\}$ stands for the diagonal (or block diagonal) matrix 
with the entries $d_1, ...$ on the main diagonal. By $\|A\|$ and by $\s(A)$, we denote the $l^2$-induced norm and  the spectrum, respectively, of some matrix $A$. We say that the matrix $X$ is positive (positive definite)
and write $X>0$ if $X$ is Hermitian (i.e., $X=X^*$) and all the eigenvalues of $X$ are positive.

\section{Preliminaries}\label{Prel}
\setcounter{equation}{0}
\subsection{Rational functions}
Recall that  a rational matrix function is called {\it strictly proper} if it tends to zero
at infinity. It is well-known \cite{KFA, LR} that  such an $m_2
\times m_1$ matrix function $\vp$ can be represented in the form
\begin{equation}
\label{app.1} \vp(z)=\clc(z I_n- \cla)^{-1}\clb,
\end{equation}
where $\cla$ is a square matrix of some order $n$, and the matrices $\clb$
and $\clc$ are of sizes $n \times m_1$ and $m_2 \times n$,
respectively. The representation (\ref{app.1})
is called a {\it realization}  of $\vp$, and  the realization
(\ref{app.1}) is said to be {\it minimal} if 
 $n$ is minimal among all possible realizations of $\vp$.
This minimal $n$  is called the {\it McMillan degree} of $\vp$. 
The
realization (\ref{app.1}) of $\vp$ is minimal if and only if
\begin{equation}
\label{app.2} {\mathrm{span}}\bigcup_{k=0}^{n-1}\im
\cla^k\clb=\BC^n,\quad {\mathrm{span}}\bigcup_{k=0}^{n-1}\im  (\cla^*)^k\clc^*=\BC^n, \quad
n= \mathrm{ord}(\cla),
\end{equation}
where Im stands for image and $\mathrm{ord}(\cla)$ stands for the order of $\cla$.
If for a pair of matrices $\{\cla, \, \clb\}$ the first equality in
(\ref{app.2}) holds, then the pair $\{\cla, \, \clb\}$ is called {\it
controllable}. If the second equality in
(\ref{app.2}) is fulfilled, then the pair $\{\clc, \, \cla\}$ is said to be {\it
observable}. 

Now, let matrix functions $\vp$ be contractive, that is, let $\vp^*\vp \leq I_{m_1}$ (or, equivalently, $\vp \vp^* \leq I_{m_2}$) hold. From \cite[Theorems 21.1.3, 21.2.1]{LR} (see also \cite[p. 191]{FKRS2013}),   the next proposition
easily follows.
\begin{Pn}\label{PnRic1} Assume that  $\vp$ is a strictly proper rational matrix function, which is contractive on $\BR$ and has no poles
in $\BC_+$, and let realization \eqref{app.1} be its minimal realization.

Then, there is a positive solution $X>0$ of the Riccati equation 
\begin{align} \label{app.3}&
 X\clb \clb^*X+\I(\cla^*X-X\cla )+\clc^*\clc =0.
\end{align}
\end{Pn}
Clearly, under conditions of Proposition \ref{PnRic1}, $\vp(z)$ is contractive on $\BC_+\cup \BR$.

In the case of the skew-self-adjoint Dirac system, we obtain Ricatti equation with minus before $\clb \clb^*$:
\begin{align} \label{app.4}&
X\clc^*\clc X+\I(\cla  X-X\cla^*)-\clb \clb^*=0.
\end{align}
This case should be dealt with in a different way and we shall do it separately
(in the next paper).

%%%%%%%%%%%%%%%%%%%%%%%%%%%%%%%%%%%%%%%%%%
%%%%%%%%%%%%%
\subsection{System \eqref{1.1}: Weyl function and inverse problem}\label{Sub2.2}
Recall that $Y(x,z)$ is the normalized by $Y(0,z)=I_m$ fundamental solution of Dirac system \eqref{1.1}, where $j$ and $V$ have the forms \eqref{1.2}.
\begin{Dn}  \label{defWeyl1} 
An $m_2 \times m_1$ matrix function $\varphi(z)$ $(z \in \BC_+)$ such that
\begin{align}&      \label{W1}
\int_0^{\infty}
\begin{bmatrix}
I_{m_1} & \vp(z)^*
\end{bmatrix}
Y(x,z)^*Y(x,z)
\begin{bmatrix}
I_{m_1} \\ \vp(z)
\end{bmatrix}dx< \infty  
\end{align}
is called a Weyl function  of the Dirac  system \eqref{1.1}  on $[0, \, \infty)$.
\end{Dn}
\begin{Rk} \label{Rk2.3} According to \cite[Sect. 2]{FKRS2013} and\cite[Sect. 2.2]{SaSaR}, the Weyl function $\vp(z)$ of the Dirac  system \eqref{1.1} always exists and is unique.
Moreover,  $\vp(z)$  is holomorphic and contractive in $\BC_+$. 
\end{Rk}
If $\vp$ is rational,  it can be prolonged (from $\BC_+$) on $\BR$ and $\BC_-$ in a natural way.
Each potential $v$ corresponding to a strictly proper rational Weyl function is generated by a fixed value $n\in \BN$
and by a quadruple of matrices, namely, by two $n\times n$ matrices $\a$ and $S_0>0$
and by  $n\times m_k$ matrices $\vt_k$ $(k=1,2)$ such that the matrix identity
\begin{align} \label{2.7}&
\a S_0-S_0\a^*=\I(\vt_1\vt_1^*- \vt_2\vt_2^*) 
\end{align}
holds. Such potentials $v$ have the form
\begin{align} \label{2.5}&
v(x)=-2 \I \vt _{1}^{\, *}\E^{\I x \alpha^{*}} S
(x)^{-1}\E^{\I x
\alpha }\vt_2,
\\ \label{2.6!}& 
S(x)=S_0+ \int_{0}^{x} \Lambda(t) \Lambda (t)^{*}dt,
\quad
\Lambda (x)= \begin{bmatrix}  \E^{- \I x \alpha } \vt_{1} 
& \E^{\I x
\alpha } \vt_2 \end{bmatrix}.
\end{align}
\begin{Dn} \label{DnPE}\cite{GKS1, FKRS2013} The potentials $v$ generated by the quadruples $\quad$ $\{\a, S_0, \vt_1, \vt_2\}$
$($where $S_0>0$ and \eqref{2.7} holds$)$
via equalities \eqref{2.5} and \eqref{2.6!},
are called 
pseudo-exponential.
\end{Dn}
\begin{Tm}\label{TmDp}\cite{FKRS2013} Let Dirac system with a pseudo-exponential
potential $v$ be given on $[0,\, \infty)$ and let $v$ be generated by the quadruple $\{\a, S_0, \vt_1, \vt_2\}$.
Then the Weyl function $\vp$ of this system has the form
\begin{align} \label{2.5!}&
\vp(z)=-\I\vt_2^*S_0^{-1}(zI_n-\t)^{-1}\vt_1, \quad \t=\a-\I\vt_1\vt_2^*S_0^{-1}.
\end{align}
\end{Tm}
The following theorem (i.e., \cite[Theorem 3.4]{FKRS2013}) presents a procedure of explicit solution
of the inverse problem (see also \cite[Theorem 5.2]{GKS6} for the $m_1=m_2$ case), which is  basic for this paper.

\begin{Tm}\label{TmIpes} Let $\vp(z)$ 
be a strictly proper rational matrix function, which  is contractive on $\BR$ and has no poles in $\BC_+$.  Assume that \eqref{app.1} is its minimal realization
and that $X>0$ is a solution of \eqref{app.3}. 

Then $\vp(z)$ is the Weyl function of the Dirac system \eqref{1.1}, the potential $v$ of which
has the form \eqref{2.5}, \eqref{2.6!}, where the quadruple  $\{\a, S_0, \vt_1, \vt_2\}$ is given $($in terms of $\cla$, $\clb$, $\clc$ and $X)$ by the relations
\begin{align} 
\label{2.6}&
\a=\cla +\I \clb \clb^*X,  \quad S_0=X^{-1}, \quad \vt_1=\clb, \quad \vt_2=-\I X^{-1} \clc^*.
\end{align}
\end{Tm}
In particular, the identity \eqref{2.7} easily follows from \eqref{app.3} and \eqref{2.6}. The uniqueness of our explicit solution of the inverse
problem is immediate from a much more general uniqueness result.

\begin{Pn} \label{PnUniq1} \cite{ALSJSp} The solution of the inverse problem to recover system \eqref{1.1} from its Weyl function
is unique in the class of Dirac systems
with the locally square summable potentials.
\end{Pn}

\begin{Rk}\label{contra}
We note that there are many quadruples generating the same pseudo-exponential potential. The quadruples,
which are recovered using \eqref{2.6}, have an important additional property$:$ controllability of the pair
$\{\a, \vt_1\}$. This property is immediate from the controllability of the pair
$\{\cla, \clb\}$.
\end{Rk}
Furthermore, the matrices $\cla$, $\clb$ and $\clc$ in the minimal realizations \eqref{app.1} of  $\vp$ are unique up to {\it basis} (similarity)
transformations:
\begin{align} \label{2.8}&
\wh{\cla}=\clt^{-1}\cla \clt, \quad \wh{\clc}=\clc \clt, \quad \wh{\clb}=\clt^{-1}\clb,
\end{align}
where $\clt$ are invertible $m\times m$ matrices. Choosing the realization of $\vp$
with $\wh \cla$, $\wh \clb$ and $\wh \clc$ instead of $\cla$, $\clb$ and $\clc$, and
adding the sign "$\, \wh{}\, $" in the notations of the corresponding matrices $\a$, $\vt_i$ and $X$, we derive
\begin{align} \label{2.9}&
\wh{\a}=\clt^{-1}\a \clt, \quad \wh{\vt}_i=\clt^{-1}\vt_i \quad (i=1,2), \quad \wh{X}=\clt^* X \clt .
\end{align}
Setting $\clt= X^{-1/2}U^*$,  where $U$ is unitary, we have $\wh X=I_m$. Hence, \eqref{2.7}
takes the form $\wh \a - \wh \a^*=\I(\wh \vt_1 \wh \vt_1^*- \wh \vt_2 \wh \vt_2^*)$.
Moreover, for the case $m_1=m_2=p$, it was shown in \cite{GKS4} that $U$ may be chosen in such a way that we have the block representations:
\begin{align} \label{2.10}&
\wh{\b}:=\wh \a - \I\wh \vt_1(\wh \vt_1 +\wh \vt_2)^*=\begin{bmatrix} \wt \b & 0 \\ 0 & \zeta \end{bmatrix}, \quad
\wh \vt_1=\begin{bmatrix} \wt \vt_1  \\  \om \end{bmatrix}, \quad \wh \vt_2=\begin{bmatrix} \wt \vt_2  \\ - \om \end{bmatrix},
\end{align}
where 
\begin{align} \label{2.11}&
\zeta=\zeta^*=\diag\{t_1 I_{n_1}, t_2 I_{n_2}, \ldots,  t_k I_{n_k}\},  \quad \s(\wt \b)\in \BC_-,
\end{align}
and $\om$ is some $\wt n\times p$ matrix, $\wt n:=\sum_{i=1}^k n_k$.

Now, introduce Dirac operator $H$ associated with the differential expression
\begin{align} &       \label{2.12}
\big(\clh f\big)(x)=\left(-\I j \frac{d}{d x}-V(x)\right)f(x),
\end{align} 
the domain of which consists of all absolutely continuous functions  $f$ from  $L^2_m(0,\infty)$, such that $\clh f\in  L^2_m(0,\infty)$
and the initial condition 
$$\begin{bmatrix} I_p  &  - I_p\end{bmatrix}f(0)=~0$$ 
holds.  Using \eqref{2.10}, it is shown in \cite{GKS4} (see also \cite[Sect. 2]{GKS6})
that the real eigenvalues of  $H$ are concentrated at the points $t_k$ and have multiplicities $n_k$, whereas the continuous 
spectrum of $H$ is described by $\wt \b$, $\wt \vt_1$ and $\wt \vt_2$. Namely, the spectral density $\vr$ of $H$ has the form
\begin{align} &       \label{2.13}
\vr(t)=g(t)^*g(t), \quad g(t):=I_p-\I(\wt \vt_1 +\wt \vt_2)^*(tI_{n-\wt n}-\wt \b)^{-1}\wt \vt_1.
\end{align} 
In view of the mentioned above connections between the quadruple $\quad$ $\,$ $\{\a, S_0, \vt_1, \vt_2\}$ and the corresponding Weyl and spectral functions,
we can consider this quadruple as some  spectral data.
%%%%%%%%%%%%%%%%%%%%%%%%%%%%%%%%%%%%%%%%%%%%%%%%%%%%%%%%%%%%%%%
%%%%%%%%%%%%%%%%%%%%%%%%%%%%%%%%%%%%%%%%%%%%%%%%%%%%%%%%%%%%%%%
%%%%%%%%%%%%%%%%%%%%%%%%%%%
\section{Stability of explicit solutions}\label{Sta}
\setcounter{equation}{0}
{\bf 1.} The following lemma is a stronger statement than Proposition \ref{PnRic1}. (Theorem 7.4.2 from \cite{LR} is used for its proof in addition to 
the Theorems 21.1.3 and 21.2.1 from
\cite{LR}, which yield Proposition \ref{PnRic1}.)
\begin{La}\label{LaJMAA} \cite{FKSjmaa}
Assume that a strictly proper rational $m_2 \times m_1$ matrix function $\vp(z)$ is contractive
on $\BR$, and that \eqref{app.1} is its minimal realization.

Then there is a unique Hermitian
solution $X$ of the Riccati equation \eqref{app.3} such that 
the relation
\begin{align} &       \label{3.1}
\s(\cla+\I \clb\clb^* X)\subset (\BC_-\cup \BR)
\end{align} 
holds. This solution $X$ is always
invertible. It is also positive if and only if  $\vp(z)$ is contractive in $\BC_+$.
\end{La}
Further, in our procedure to recover the potential $v$ from the Weyl function $\vp$, we shall look for this
particular solution $X$ of \eqref{app.3}. More precisely, we start with  the strictly proper rational $m_2 \times m_1$ matrix function $\vp(z)$,
which is contractive on $\BR$ and has no poles in $\BC_+$. Hence, $\vp(z)$ is contractive in $\BC_+$, and so, according to Lemma \ref{LaJMAA},
we have $X>0$.
By $\clg_n$ we denote the class of triples $\{\wt \cla, \wt \clb, \wt \clc\}$ which determine minimal realizations $\wt \vp(z)=\wt \clc(zI_n-\wt \cla)^{-1}\wt \clb$
of $m_2\times m_1$ matrix functions $\wt \vp(z)$ contractive on $\BR\cup \BC_+$. First, we consider the stability of recovery of $X$ from
$\{ \cla,  \clb,  \clc\}\in \clg_n$.
\begin{Dn}\label{stabX} The recovery of $X>0$, satisfying \eqref{3.1},  from the minimal realization \eqref{app.1} $($where $\{ \cla,  \clb,  \clc\}\in \clg_n)$  of $\vp(z)$ is called  stable
if for any $\ve>0$ there is $\delta >0$ such that for each $\{\wt \cla, \wt \clb, \wt \clc\}$ satisfying conditions
\begin{align} &       \label{3.1!}
\{\wt \cla, \wt \clb, \wt \clc\}\in \clg_n, \quad \|\cla -\wt \cla\|+ \| \clb -\wt \clb\|+\|\clc- \wt \clc\|<\delta
\end{align} 
there is a solution $\wt X=\wt X^*$ of the equation $ \wt X \wt \clb \wt \clb^* \wt X+\I(\wt \cla^* \wt X-\wt X\wt \cla )+\wt \clc^* \wt \clc =0$
in the neighborhood $\| X-\wt X\|<\ve$ of $X$.
\end{Dn}
The stability of solutions $X$ of an important class of Riccati equations was shown in \cite[Theorem 4.4]{RaRo} for a somewhat wider class of
perturbations than described in our definition and we shall use this theorem in order to prove our first stability statement.
\begin{Tm}\label{TmStab1} The recovery of $X>0$, satisfying \eqref{3.1},  from the minimal realization \eqref{app.1} of $\vp(z)$ $($with $\{ \cla,  \clb,  \clc\}\in \clg_n)$
is stable.
\end{Tm}
\begin{proof}.
Assuming that a minimal realization \eqref{app.1} of $\vp(z)$
is given (that is, matrices $\cla$, $\clb$ and $\clc$ are given), we consider equation \eqref{app.3}. Putting
\begin{align} &       \label{3.2}
A_0=-\I(\cla + c I_n) \quad (c \in \BR), \quad C_0=-\clc^*\clc, \quad D_0=\clb\clb^*,
\end{align} 
we see that the equation \eqref{app.3} coincides with the Riccati equation
$XD_0X+XA_0+A_0^*X-C_0=0$ considered in \cite[Subsection 4.2]{RaRo}.

Now, we deal with the conditions (i)-(iv) (on the coefficients $A_0, C_0, D_0$) from \cite[Subsection 4.2]{RaRo}.
(Only perturbations satisfying these conditions are allowed   in \cite[Subsection 4.2]{RaRo} and we will
show that  the conditions $\quad$ (i)-(iv) are satisfied in the case $\{ \cla,  \clb,  \clc\}\in \clg_n$.)
Equalities \eqref{3.2} and the fact that the pair $\{\cla, \clb\}$ is controllable
imply that conditions (i) and (ii)  in \cite[Subsection 4.2]{RaRo} are fulfilled.
In a similar way we derive that conditions (i) and (ii) are fulfilled for the Ricatti
equations $\wt X \wt D_0 \wt X+\wt X\wt A_0+\wt A_0^*\wt X-\wt C_0=0$
corresponding to all the triples $\{\wt \cla, \wt \clb, \wt \clc\}\in \clg_n$.
For sufficiently large values of $|c|$, the requirement (iii)
that the matrix
\begin{align} &       \label{3.3}
H=\begin{bmatrix}-C_0 & A_0^* \\ A_0 & D_0 \end{bmatrix}
\end{align} 
satisfies the condition $\det H\not= 0$ and that signature of $H$ equals zero is also
fulfilled. Clearly, $c$ may be chosen so that (iii) is valid if we substitute 
the triple $\{\cla, \clb, \clc\}$ with any triple $\{\wt \cla, \wt \clb, \wt \clc\}\in \clg_n$ in some small neighborhood of $\{\cla, \clb, \clc\}$.
Finally, according to Lemma \ref{LaJMAA}, there are hermitian solutions of equations $\wt X \wt D_0 \wt X+\wt X\wt A_0+\wt A_0^*\wt X-\wt C_0=0$,
that is, condition (iv) holds. Since conditions (i)-(iv)  from \cite[Subsection 4.2]{RaRo} are fulfilled, the stability with respect to perturbations in the class $\clg_n$ 
will follow from the stability in the sense of \cite[Theorem 4.4]{RaRo}.

Using again Lemma \ref{LaJMAA}, we choose the solution $X>0$ of \eqref{app.3} satisfying \eqref{3.1}. 
It is immediate that one of the
equivalent statements from \cite[Theorem 4.4]{RaRo} is valid for our $X$.
That is, according to \eqref{3.1} and \eqref{3.2}, the equality $\Im(\la)=0$ holds for each $\la$ from the set
\begin{align}       \nn
\s\big(\I(A_0&+D_0X)\big)\cap\s\big(-\I(A_0^*+XD_0)\big)
\\ & \label{3.4}
=\s\big(\cla+\I \clb\clb^*X+cI_n\big)\cap\s\big((\cla+\I \clb\clb^*X+cI_n)^*\big),
\end{align} 
and so the statement (d) from \cite[Theorem 4.4]{RaRo} holds.
Therefore, by virtue of \cite[Theorem 4.4]{RaRo}, $X$ is a stable and isolated solution of \eqref{app.3}.
\end{proof}
{\bf 2 .}  Now, we will show that small perturbations of the quadruple $\{\a, S_0, \vt_1,\vt_2\}$ result in small
perturbations of the corresponding potential $v$. We note that we consider only perturbations which do not change $m_1, m_2$ and $n$ .
\begin{Dn}\label{adm} The quadruple $\{\a, S_0, \vt_1,\vt_2\}$ is called {\rm admissible} if $S_0>0$ and \eqref{2.7} holds,
and it is called {\rm spectral} if it is admissible, the pair $\{\a, \vt_1\}$ is controllable and 
\begin{align}       & \label{3.5}
\s(\a)\subset (\BR\cup \BC_-).
\end{align} 
\end{Dn}
\begin{Rk}\label{RkSpeq} Theorem \ref{TmDp} and Remark \ref{Rk2.3} show that the Weyl function corresponding to any
pseudo-exponential potential is rational and contractive. Then Theorem \ref{TmIpes}, Proposition \ref{PnUniq1} and
Lemma \ref{LaJMAA}  imply that this potential $($uniquely recovered from the Weyl function$)$ is generated, in particular, by
a spectral quadruple. In other words, each pseudo-exponential potential is generated by some spectral quadruple.
\end{Rk}
\begin{Tm} \label{TmStab2}
Let a spectral quadruple $\{\a, S_0, \vt_1,\vt_2\}$ be given. Then, for any $\ve >0$ there is $\delta >0$ such that
each pseudo-exponential potential $\wt v$ generated by an admissible quadruple $\{\wt \a,\wt S_0,\wt \vt_1,\wt \vt_2\}$ satisfying condition
$$\|\a-\wt \a\|+\|S_0-\wt S_0\|+ \|\vt_1-\wt \vt_1\|+\|\vt_2-\wt \vt_2\|<\delta$$
belongs to the $\ve$-neighborhood of $v$ generated by $\{\a, S_0, \vt_1,\vt_2\}$, that is, 
\begin{align}       & \label{3.6}
\sup_{x\in [0,\infty)}\|v(x)-\wt v(x)\|<\ve.
\end{align} 
\end{Tm}
In order to prove the theorem above we generalize (for the case when $m_1$ does not necessarily
equal $m_2$ and $S_0$ does not necessarily equal $I_n$) some results from \cite{GKS6} on asymptotics of 
\begin{align}       & \label{3.7}
Q(x):=S_0+2 \int_0^x
\E^{2\I t \a}\vt_2\vt_2^*\E^{-2\I t \a^*}dt.
\end{align} 
\begin{La}\label{LaAsymp} The following relations are valid for a spectral quadruple\\ $ \{\a, S_0, \vt_1,\vt_2\}:$
\begin{align}       & \label{3.8}
\lim_{x\to \infty}Q(x)^{-1}=0, \quad
\lim_{x\to \infty}\|Q(x)^{-1}\E^{2\I x \a}\vt_2\|=0.
\end{align} 
\end{La}
\begin{proof}.  The proof uses some steps from the proof  of \cite[Proposition 3.3]{GKS6}.
Since $Q(x)$ is increasing and is bounded below by $S_0>0$, there is a limit $\ka_Q:=\lim_{x\to \infty}Q(x)^{-1}$.
Next, we prove that $\ka_Q=0$.  From the definition \eqref{3.7} and identity \eqref{2.7} we derive
\begin{align}       \nn
\a Q(x)-Q(x)\a^*&=\a S_0-S_0\a^*-\I\left(\E^{2\I x\a}\vt_2\vt_2^*\E^{-2\I x \a^*}-\vt_2\vt_2^*\right)
\\ & \label{3.9}
=\I \vt_1\vt_1^*-\I \E^{2\I x\a}\vt_2\vt_2^*\E^{-2\I x \a^*}.
\end{align} 
Multiplying (from both sides) the left-hand side and right-hand side of \eqref{3.9} by $Q(x)^{-1}$, we obtain
\begin{align}     &\nn
Q(x)^{-1}\a -\a^*Q(x)^{-1}-\I Q(x)^{-1}\vt_1\vt_1^*Q(x)^{-1}
\\  & \label{3.10}  \qquad
=-\I Q(x)^{-1}\E^{2\I x\a}\vt_2\vt_2^*\E^{-2\I x \a^*} Q(x)^{-1}.
\end{align} 
Passing in \eqref{3.10} to the limit, we see that
\begin{align}     
& \label{3.11}
\lim_{x\to \infty}Q(x)^{-1}\E^{2\I x\a}\vt_2\vt_2^*\E^{-2\I x \a^*} Q(x)^{-1}=\I(\ka_Q\a-\a^*\ka_Q-\I\ka_Q\vt_1\vt_1^*\ka_Q).
\end{align} 
On the other hand, formula \eqref{3.7} yields
$$\frac{d}{dx}Q(x)^{-1}=-2Q(x)^{-1}      \E^{2\I x \a}\vt_2\vt_2^*\E^{-2\I x \a^*}        Q(x)^{-1},$$
and so we have
\begin{align}     
& \label{3.12}
\int_0^{\infty}Q(t)^{-1}      \E^{2\I t \a}\vt_2\vt_2^*\E^{-2\I t \a^*}        Q(t)^{-1}dt=\frac{1}{2}(S_0^{-1}-\ka_Q)<\infty.
\end{align} 
Taking into account \eqref{3.12} and the fact that there exists a limit of the expression integrated in \eqref{3.12}
(see \eqref{3.11}), we derive that this limit equals zero. That is, we rewrite \eqref{3.11} in the form
\begin{align}     
& \label{3.13}
\ka_Q\a-\a^*\ka_Q-\I\ka_Q\vt_1\vt_1^*\ka_Q=0.
\end{align}
Moreover, since the left-hand side in \eqref{3.11} tends to zero, the second equality in \eqref{3.8} is already proved.

Recall that the first equality in \eqref{3.8} is equivalent to $\ka_Q=0$. Now, we prove $\ka_Q=0$ by negation.
For this, we rewrite \eqref{3.13} in the form $\a^*\ka_Q=\ka_Q(\a -\I \vt_1\vt_1^*\ka_Q)$, which implies that the range
of $\ka_Q$ is an invariant subspace of $\a^*$. Thus, assuming $\ka_Q\not= 0$, we obtain
that there is an eigenvector $\ka_Q g$ of $\a^*$:
$\,
\a^*\ka_Q g=c \ka_Q g, \quad \ka_Q g\not=0, \quad g\in \BC^n.
$

Finally, consider the expression $\I g^*(\ka_Q\a-\a^*\ka_Q)g$. In view of \eqref{3.5}, for the eigenvalue $c$ of $\a^*$
we have $\Im (c)\geq 0$, and so 
$$\I g^*(\ka_Q\a-\a^*\ka_Q)g=\I(\ov c- c)g^*\ka_Q g\geq 0.$$ 
On the other hand,
we have $\vt_1^*\ka_Q g\not=0$ because the pair $\{\a,\vt_1\}$ is controllable. Hence, the inequality
 $\I g^*(\ka_Q\a-\a^*\ka_Q)g=-g^*\ka_Q\vt_1\vt_1^*\ka_Qg<0$ follows from \eqref{3.13}.
 We arrive at a contradiction, that is, $\ka_Q=0$.
\end{proof}
In the case of admissible quadruples $\{\a, S_0, \vt_1,\vt_2\}$, the matrix identity 
\begin{align}     
& \label{3.14}
\a S(x)-S(x)\a^*=\I \Lam(x)j\Lam(x)^*
\end{align}
(see \cite[formula (3.6)]{FKRS2013}) coincides  with \eqref{2.7} at $x=0$ and easily follows  from \eqref{2.7} and \eqref{2.6!} for $x>0$.
In other words, $\a$, $S(x)$ and $\Lam(x)$ form an $S$-node (and, moreover, the so called Darboux matrix function corresponding
to $v(x)$ coincides with the transfer matrix function \cite{SaL1, SaL2, SaSaR} in Lev Sakhnovich sense).
Using \eqref{2.6!}, \eqref{3.7} and \eqref{3.14}, we derive
\begin{align}     
& \nn
Q^{\prime}(x)=\big(\E^{\I x \a}S(x)\E^{-\I x \a^*}\big)^{\prime}, \quad Q(0)=S(0) \quad \left(Q^{\prime}:=\frac{d}{dx}Q\right),
\end{align}
and so the following equality is valid:
\begin{align}     
& \label{3.15}
Q(x)=\E^{\I x \a}S(x)\E^{-\I x \a^*}.
\end{align}
{\it Proof of Theorem \ref{TmStab2}.} Now, we consider a pseudo-exponential potential $v$ generated by the spectral quadruple 
$\{\a, S_0, \vt_1,\vt_2\}$ and pseudo-exponential potentials $\wt v$ generated by admissible quadruples $\{\wt \a,\wt S_0,\wt \vt_1,\wt \vt_2\}$
belonging to a neighborhood of  $\{\a, S_0, \vt_1,\vt_2\}$. The matrix function $Q$ corresponding to $\{\wt \a,\wt S_0,\wt \vt_1,\wt \vt_2\}$ is
denoted by $\wt Q$. In view of \eqref{2.5} and \eqref{3.15}, we have:
\begin{align}     
& \label{3.16}
v(x)=-2\I \vt_1^*Q(x)^{-1}\E^{2 \I x \a}\vt_2, \quad \wt v(x)=-2\I \wt \vt_1^*\wt Q(x)^{-1}\E^{2 \I x \wt \a}\wt \vt_2.
\end{align}
It is immediate from the proof of Lemma \ref{LaAsymp} that \eqref{3.10} holds for admissible quadruples as well. That is, 
we may rewrite \eqref{3.10} for $\{\wt \a,\wt S_0,\wt \vt_1,\wt \vt_2\}$:
\begin{align}     &\nn
\wt Q(x)^{-1}\wt \a -\wt \a^*\wt Q(x)^{-1}-\I \wt Q(x)^{-1}\wt \vt_1 \wt \vt_1^* \wt Q(x)^{-1}
\\  & \label{3.17}  \qquad
=-\I \wt Q(x)^{-1}\E^{2\I x\wt \a}\wt \vt_2 \wt \vt_2^*\E^{-2\I x\wt \a^*}\wt Q(x)^{-1}.
\end{align} 
Since $Q$ and $\wt Q$ are monotonic and the first relation in \eqref{3.8} is valid, we may choose $x_0>0$ and some neighborhood
of  $\{\a, S_0, \vt_1,\vt_2\}$ so that $Q(x)$ and $\wt Q(x)$ are large enough for $x\geq x_0$, and so the left-hand sides of 
\eqref{3.10} and \eqref{3.17} are small enough. Hence,  the right-hand sides of 
\eqref{3.10} and \eqref{3.17}  are also small enough. Therefore, taking into account \eqref{3.16}, we see that for any $\ve>0$
there are $x_0>0$ and $\delta_1>0$ such that the next inequality holds in the $\delta_1$-neighborhood 
$(\|\a-\wt \a\|+\|S_0-\wt S_0\|+ \|\vt_1-\wt \vt_1\|+\|\vt_2-\wt \vt_2\|<\delta_1)$ of  $\{\a, S_0, \vt_1,\vt_2\}$:
\begin{align}       & \label{3.18}
\sup_{x\in [x_0,\infty)}\|v(x)-\wt v(x)\|<\ve .
\end{align} 
It easily follows from the definition of $Q$ and $\wt Q$ and from \eqref{3.16} that there is some $\delta_2$-neighborhood
of $\{\a, S_0, \vt_1,\vt_2\}$, where we have
\begin{align}       & \label{3.19}
\sup_{x\in [0,x_0)}\|v(x)-\wt v(x)\|<\ve .
\end{align} 
Clearly, inequalities \eqref{3.18} and \eqref{3.19} yield \eqref{3.6}  (for $\delta=\min(\delta_1,\delta_2)$). $\blacksquare$

\begin{Rk} It is immediate from the second relation in \eqref{3.8}, the first relation in \eqref{3.16} and Remark \ref{RkSpeq} that all pseudo-exponential
potentials tend to zero at infinity.
\end{Rk}
{\bf 3.}  Theorems \ref{TmIpes}, \ref{TmStab1} and \ref{TmStab2} as well as Lemma \ref{LaJMAA} yield the result below on the stability of the procedure of 
solving inverse problem.

\begin{Tm} The procedure $($given in Theorem \ref{TmIpes}$)$ to uniquely recover the pseudo-exponential potential $v$ of Dirac system \eqref{1.1} from a minimal
realization of the Weyl function
$($i.e., of some strictly proper rational $m_2\times m_1$  matrix function, which is contractive in $\BC_+)$ is stable once we agree to choose such a positive solution $X$
of the Riccati equation \eqref{app.3}  that \eqref{3.1} holds.
\end{Tm} 

%%%%%%%%%%%%%%%%%%%%%%%%%%%%%%%%%%%%%%%%%%%%%%%%%%%%%%%%%%%%%%%%%%
%%%%%%%%%%%%%%%%%%%%%%%%%%%%%%%%%%%%%%%%%%%%%%%%%%%%%%%%%%
%%%%%%%%%%%%%%%%%%%%%%%%%%%%%%%%%%

\bigskip 
\noindent{\bf Acknowledgments.}
 {This research  was supported by the
Austrian Science Fund (FWF) under Grant  No. P24301.}
%%%%%%%%%%%%%%%%%%%%%%%%%%%%%%%%%%%%%%%%%%%%%%
%%%%%%%%%%%%%%%%%%%%%%%%%%%%%%%%%%%%%%%%%%%%%%%
\newpage
%%%%%%%%%%%%%%%%%%%%%%%%%%%%%%%%%%%%%%%%%%%
%%%%%%%%%%%%%%%%%%%%%%%%%%%%%%%%%%%%%%%%%%%%

\begin{flushright}
Alexander Sakhnovich,\\
Vienna
University
of
Technology,
Austria, \\
e-mail: {\tt oleksandr.sakhnovych@tuwien.ac.at}

\end{flushright}

%%%%%%%%%%%%%%%%%%%%%%%%

\end{document}